\documentclass[superscriptaddress,twocolumn,aps,jcp,amsmath,amssymb,reprint]{revtex4-1}

\usepackage{graphicx}
\usepackage{dcolumn}
\usepackage{bm}

\usepackage{xcolor}

\DeclareMathOperator{\divop}{div}
\DeclareMathOperator{\grad}{grad}

\newcommand{\ud}{\,\mathrm{d}}
\newcommand{\cP}{\mathcal{P}}

\newcommand{\RR}{\mathbb{R}}

\newcommand{\PP}{\mathbb{P}}

\newcommand{\bd}[1]{\boldsymbol{#1}}

\newcommand{\abs}[1]{\lvert#1\rvert}

\newcommand{\average}[1]{\langle#1\rangle}

\newcommand{\barint}{\kern4pt \raise3.4pt\hbox{\vrule height.6pt
		width7pt} \kern-11pt \int}

\newcommand{\figref}[1]{Fig.~\ref{#1}}

\begin{document}


\title{Integrated Tempering Enhanced Sampling Method as the Infinite Switching Limit of Simulated Tempering}%

\author{Zhiyi You}%
\affiliation{Department of Statistics, University of California,
  Berkeley, CA 94720, United States }%

\author{Liying Li}%
\affiliation{Courant Institute of Mathematical Sciences, New York
  University, New York, N.Y. 10012, United States }%

\author{Jianfeng Lu}%
\email{jianfeng@math.duke.edu} \affiliation{ Department of
  Mathematics, Department of Physics, and Department of Chemistry,
  Duke University, Durham, NC 27708, United States }%

\author{Hao Ge}%
\email{haoge@pku.edu.cn} \affiliation{Beijing International Center for
  Mathematical Research (BICMR) and Biodynamic Optical Imaging Center (BIOPIC), Peking
  University, Beijing 100871, P.R.China }%

\date{\today}

\begin{abstract}
  Fast and accurate sampling method is in high demand, in order to
  bridge the large gaps between molecular dynamic simulations and
  experimental observations. Recently, integrated tempering enhanced
  sampling method (ITS) has been proposed and successfully applied to
  various biophysical examples, significantly accelerating
  conformational sampling. The mathematical validation for its
  effectiveness has not been elucidated yet. Here we show that the
  integrated tempering enhanced sampling method can be viewed as a
  reformulation of the infinite switching limit of simulated tempering
  method over a mixed potential. Moreover, we demonstrate that the
  efficiency of simulated tempering molecular dynamics (STMD) improves
  as the frequency of switching between the temperatures is increased,
  based on the large deviation principle of empirical
  distributions. Our theory provides the theoretical justification of
  the advantage of ITS.  Finally, we illustrate the utility of the
  infinite switching simulated tempering method through several
  numerical examples.

\end{abstract}

\maketitle


\section{Introduction}

Molecular dynamics simulation is a powerful method for investigating
microscopic biochemical systems. However, when the system contains
barriers due to energy or entropy, direct molecular dynamics
simulation will lead to a very high computational cost. Therefore,
people turn to use the enhanced sampling method, sacrificing the real
dynamic information in order to capture the desired Gibbs
distribution. Until now, sampling methods have already been
successfully applied in multiple disciplines including statistical
physics, chemical physics, Bayesian statistics, machine learning, and
related areas.


The efficiency of the sampling approach highly depends on the convergence rate to the thermodynamic equilibrium. 
Tempering schemes, such as simulated tempering molecular dynamics
(STMD) \cite{marinari1992simulated, kerler1994simulated} and replica
exchange molecular dynamics (REMD) \cite{hansmann1997parallel,
  sugita1999replica, earl2005parallel}, were developed and are among the most
popular methods to overcome the metastability and to enhance the
convergence to equilibrium, thanks to their effectiveness and simplicity to use. The
basic idea of both is to use one or more artificial high temperatures
to accelerate exploration of the conformational space, while using the
physical temperature to sample the desired physical observable. The
interaction between different temperatures is designed so to guarantee
the unbiasedness of the estimation. These two methods are quite comparable to each other. STMD gives a higher rate of delivering the system between
different temperature states as well a higher rate of traversing the
energy space \cite{zhang2008comparison}, but requires the estimation of partition function first.
  
Another tempering algorithm named ``integrated tempering enhanced
sampling'' (ITS), which was recently introduced by Gao
\cite{gao2008integrate}, uses a temperature-biased effective potential
energy to run the MD simulation and has been successfully applied to
various biochemical examples \cite{YangShaoGao:09,YangLiu:2015}.
Noting that ITS also uses auxiliary temperatures to accelerate the
convergence, it is desirable to compare it with previous methods, and
find out why it performs better in many cases, which has not been
studied theoretically to the best of our knowledge.

As one of the main results in this paper, we discover that the
integrated tempering enhance sampling method is in fact a reformulated
version of the infinite switching limit of STMD, that is, the limit of
the simulated tempering when the attempt switching frequency goes to
infinity. Moreover, we show that as the switching rate in STMD
increases, the empirical measure converges faster towards stationary
distribution, using the large deviation principle; and thus the
sampling efficiency of STMD increases as the switching frequency
increases. Combining the two theoretical findings, we justify the
efficiency of ITS over conventional STMD.  Finally, we compare the
infinite switching STMD against the normal STMD with two numerical
examples: an artificial high dimensional system and a more realistic
Lennard-Jones example \cite{dellago1999calculation} with 16 atoms. The
numerical results validate our theoretical findings.

Our study of the switching rate of the STMD is closely related to the
recent progress in understanding the swapping frequency of tempering
schemes, started with replica exchange methods. For REMD, it has been
discovered through numerical examples that the sampling efficiency
increases when the frequency of swapping of temperatures is pushed up
to infinity \cite{sindhikara2008exchange}; but directly increasing the
swapping rate to reach this limit is computationally infeasible, as
many swaps are needed per MD step. A breakthrough was made in
\cite{plattner2011, dupuis2012infinite} which proposed an explicit way
to reach this limit and also proved that indeed the sampling
efficiency increases in such limit. A natural reformulation of
infinite swapping REMD was later proposed in \cite{lu2013infinite}
which leads to an easy implementation as a simple patch to
conventional molecular dynamics. For multiple temperatures, an
efficient implementation based on ideas of multiscale integrator was
also proposed for infinite swapping REMD \cite{Yu2016,
  LuVandenEijnden2017}, which leads to practical applications in
sampling configurational space of large biomolecules. We note that
increasing the switching rate in the simulated tempering has been
mentioned in \cite{LuVandenEijnden2017} and
\cite{MartinssonLuLeimkuhlerVandenEijnden} under a general framework
of viewing tempering schemes as MD process augmented by jumping
processes, without much details in particular the connection with ITS.

In the remaining article, we will first recall the simulated tempering molecular dynamics as a stochastic process in Sec.~\ref{sec:STMD}. We will justify using a large switching frequency via large deviation principle for the empirical distribution in Sec.~\ref{sec:LDT} with detailed derivation given in the Appendix. Infinite switching limit of STMD and its reformulation are derived and generalized in Sec.~\ref{sec:inf} and~\ref{sec:gen}. The identification of ITS and infinite switching STMD is given and discussed in Sec.~\ref{sec:ITS}, followed by numerical examples in Sec.~\ref{sec:numerics}. Some conclusive remarks are given in Sec.~\ref{sec:conclude}.

\section{Simulated tempering molecular dynamics}\label{sec:STMD}

We begin by recalling the simulated tempering molecular dynamics.  For
the sake of simplicity we will first discuss the algorithm for system
governed by overdamped Langevin equations (i.e., when inertia can be
neglected). More general dynamics will be discussed in 
Sec.~\ref{sec:gen}. Consider 
\begin{equation}
\label{eq:overdamp}
\dot{\bd{x}} = \bd{f}(\bd{x})
+ \sqrt{2  \beta^{-1}}\, \bd{\eta},
\end{equation}
where $\bd{x} \in \RR^{3n}$ denotes the configuration of the system
with $n$ particles, $\bd{f}(\bd{x}) = -\nabla V(\bd{x})$ is the force
associated with the potential $V(\bd{x})$, $\beta= 1/k_B T$ is the
inverse temperature, and $\bd{\eta}$ is a $3n$-dimensional white-noise
with independent components. We choose a unit system so that the
friction coefficient becomes one for notational
simplicity. Eq.~\eqref{eq:overdamp} is consistent with the Boltzmann
equilibrium probability density
\begin{equation}
\label{eq:eqpdf}
\rho_\beta(\bd{x}) = Z_\beta^{-1} e^{-\beta V(\bd{x})},
\end{equation}
where the normalization constant
$Z_\beta = \int_{\RR^{3n} }e^{-\beta V(\bd{x})} d\bd{x}$ is the
partition function.

Assuming we only take two temperatures in simulated tempering (the
generalization to more temperatures will be considered in 
Sec.~\ref{sec:gen}), we
replace~\eqref{eq:overdamp} by
\begin{equation}
\label{eq:overdampmixt}
\dot{\bd{x}} = \bd{f}(\bd{x}) + \sqrt{2 \beta^{-1}(t)}\, \bd{\eta},
\end{equation}
where the inverse temperature $\beta(t)$ now attempts switches between
the physical and the artificial temperatures with frequency $\nu$, and
the attempted switch from $\beta_1$ to $\beta_2$ are accepted with
probability
\begin{equation}
\label{eq:acceptanceprob}
g_{\beta_1\beta_2}(\bd{x}) = \min\left(\frac{n_2\, e^{-\beta_2 V(\bd{x})}}{n_1\, e^{-\beta_1 V(\bd{x})}}, 1\right),
\end{equation}
where $n_1$ and $n_2$ are some weighting parameters that will be
further discussed below.

Note that in the simulated tempering overdamped dynamics, both the
configuration $\bd{x}$ and the inverse temperature $\beta$ are
dynamical variables: $\bd{x}$ follows \eqref{eq:overdampmixt},
$\beta$ follows a jump process and the two dynamics are coupled. To
understand the dynamics, it is useful to consider its infinitesimal
generator, given by (we use $\beta'$ to denote the temperature other
than $\beta$, namely $\beta' = \beta_2$ if $\beta = \beta_1$ and vice
versa)
\begin{multline}
\label{eq:infgen}
(\mathcal{L}^\nu u)(\bd{x},\beta)=
-\nabla V(\bd{x})\cdot\nabla_{\bd{x}} u(\bd{x},\beta) + \beta^{-1} \Delta_{\bd{x}} u(\bd{x},\beta)\\
-\nu g_{\beta,\beta'}(\bd{x})u(\bd{x}, \beta) +\nu
g_{\beta',\beta}(\bd{x})u(\bd{x}, \beta')
\end{multline}
for a smooth function $u:\RR^{3n}\times\{\beta_1,\beta_2\} \to \RR$. 
This means that the evolution of a physical observable $u$ under the dynamics is given by 
\begin{multline}
  \mathbb{E} \bigl( u(\bd{x}(t+\delta t), \beta(t + \delta t)) \mid \bd{x}(t), \beta(t) \bigr) = u(\bd{x}(t), \beta(t)) \\
  + \delta t (\mathcal{L}^\nu u)(\bd{x}(t),\beta(t)) + o(\delta t), 
\end{multline}
or equivalently, the density $\rho$ of  $(\bd{x}, \beta)$ evolves
under the adjoint operator:
\begin{equation}\label{eq:forward}
  \begin{aligned}
    \partial_t \rho(t, \bd{x}, \beta) & = \bigl(\mathcal{L}^\nu\bigr)^{\dagger}  \rho(t, \bd{x}, \beta)\\
    & = \Delta V(\bd{x}) \cdot \rho(t, \bd{x}, \beta)
    + \nabla V(\bd{x}) \cdot \nabla_{\bd{x}} \rho(t, \bd{x}, \beta) \\
    & \qquad + \beta^{-1} \Delta_{\bd{x}} \rho(t, \bd{x}, \beta)
    -\nu g_{\beta,\beta'} \rho(t, \bd{x}, \beta) \\
    & \qquad +\nu g_{\beta',\beta} \rho(t, \bd{x}, \beta').
    \end{aligned}
\end{equation}
For the infinitesimal generator, the first two terms on the right hand
side of \eqref{eq:infgen} correspond to the overdamped dynamics
\eqref{eq:overdampmixt} while the last two terms correspond to the
jump process of $\beta$, $-\nu g_{\beta, \beta'}(\bd{x})$ being the
rate of switching from temperature $\beta$ to $\beta'$, and
$\nu g_{\beta', \beta}(\bd{x})$ is the rate of switching to $\beta$
from the other temperature $\beta'$ (attempt switching frequency
adjusted by the acceptance rate).

It is straightforward to check the equilibrium distribution of the coupled dynamics of $(\bd{x}, \beta)$, as the stationary solution of  \eqref{eq:forward} is given by 
\begin{equation}
\label{eq:9}
{\varrho}(\bd{x},\beta) = \frac{n_1 e^{-\beta_1 V(\bd{x})} \delta_{\beta, \beta_1}
+ n_2 e^{-\beta_2 V(\bd{x})} \delta_{\beta, \beta_2}}{n_1 Z_{\beta_1} + n_2 Z_{\beta_2}},
\end{equation}
where $\delta_{\beta,\beta_1}$ denotes the Kronecker delta:
$\delta_{\beta, \beta_1} = 1$ if and only if $\beta = \beta_1$. This
is nothing but a weighted average of Boltzmann densities at
temperatures $\beta_1$ and $\beta_2$, with weight $n_1 Z_{\beta_1}$
and $n_2 Z_{\beta_2}$ respectively. On each fixed temperature, the
distribution is generated by the ordinary molecular dynamics and thus
follows Boltzmann distribution; meanwhile, at each fixed
configuration, the proportion of the two temperatures is 
given by
\begin{equation}
\label{eq:weight}
\omega_j(\bd{x}) = \frac{n_j e^{-\beta_j V(\bd{x})}} {n_1 e^{-\beta_1 V(\bd{x})} + n_2 e^{-\beta_2 V(\bd{x})}},
\quad j = 1, 2
\end{equation}
as determined by the acceptance probability
\eqref{eq:acceptanceprob}. Therefore, \eqref{eq:overdampmixt} together
with the jumping process of $\beta$ samples the
$\varrho(\bd{x}, \beta)$ as the equilibrium distribution.

Summing $\varrho(\bd{x}, \beta)$ over $\beta$ then gives the
marginal equilibrium density for the configuration position alone
\begin{equation}
\label{eq:eqpdfmixt}
\varrho(\bd{x}) = \varrho(\bd{x}, \beta_1) +\varrho(\bd{x}, \beta_2)　=　\frac{n_1 e^{-\beta_1 V(\bd{x})}
  + n_2 e^{-\beta_2 V(\bd{x})}}{n_1 Z_{\beta_1} + n_2 Z_{\beta_2}},
\end{equation}
which is a weighted average of the Boltzmann densities at the two
temperatures $\beta_1$ and $\beta_2$. As a result the ensemble average
at the physical temperature ($\beta_1$) of any observable $A(\bd{x})$
can be estimated from 
\begin{equation}
\label{eq:1eavg}
\begin{aligned}
  \average{A}_{\beta_1} &\equiv \int_{\RR^{3n}} A(\bd{x}) \rho_{\beta_1}(\bd{x}) \ud\bd{x}\\
  & = \frac{n_1 Z_{\beta_1} + n_2 Z_{\beta_2}}{n_1 Z_{\beta_1}} \int_{\RR^{3n}} A(\bd{x})\, \omega_1(\bd{x}) \varrho(\bd{x}) \ud\bd{x} \\
  & = \biggl( \int_{\RR^{3n}} \omega_1(\bd{x}) \varrho(\bd{x}) \ud \bd{x}
  \biggr)^{-1}
  \int_{\RR^{3n}} A(\bd{x})\, \omega_1(\bd{x}) \varrho(\bd{x}) \ud\bd{x} \\
  & \approx \dfrac{\displaystyle \lim_{T\to\infty} \frac1T\int_0^T A(\bd{x}(t))
    \, \omega_1(\bd{x}(t)) \ud t}{\displaystyle \lim_{T\to\infty}
    \frac1T\int_0^T \omega_1(\bd{x}) \ud t}\, , 
\end{aligned}
\end{equation}
%
where we have used that the ensemble average equals to the time
average thanks to ergodicity. Note that in principle we shall use two
independent realization to estimate the numerator and denominator on
the right hand side of \eqref{eq:1eavg} to get an unbiased estimator,
while in practice the bias is well-controlled and the variance of the
estimator usually dominates the error.

Note that the above holds for arbitrary positive weighting factors
$n_1$ and $n_2$; while the sampling efficiency of the method depends
on the choice. The conventional wisdom in simulated tempering method
is to choose $n_1$ and $n_2$ according to the partition function:
\begin{equation}
  n_1 = Z_{\beta_1}^{-1}\qquad \text{and} \qquad n_2 = Z_{\beta_2}^{-1}. 
\end{equation}
With this choice, we have 
\begin{equation}
  \varrho(\bd{x}, \beta) = \frac{1}{2} \rho_{\beta_1}(\bd{x}) \delta_{\beta, \beta_1} + \frac{1}{2} \rho_{\beta_2}(\bd{x}) \delta_{\beta, \beta_2}
\end{equation}
and thus the sampler spends equal amount time at the two
temperatures. It is also possible to choose different $n_i$ to
emphasize one of the temperatures. Of course, the partition functions
$Z_{\beta_i}$ are not known \emph{a priori}, and hence the usual
approach is to start with some initial guess of the weighting factors
and then adaptively adjust them on-the-fly using an iterative method,
see \textit{e.g.}, \cite{ParkPande:2007, gao2008integrate, Tan:17, MartinssonLuLeimkuhlerVandenEijnden}. 
 In what follows, we will mainly
focus on the choice of the switching frequency $\nu$, and thus for the
simplicity of discussion, we will assume that $n_i$ are fixed during
the sampling process and only make some remarks about adjusting $n_i$
afterwards.

\section{Large deviation principle for empirical measure}\label{sec:LDT}

As we have discussed above, for any choice of the switching frequency
$\nu$, the trajectory $(\bd x^\nu , \beta^\nu)$ (we put superscript
$\nu$ to emphasize the $\nu$ dependence) of the simulated tempering
overdamped dynamics can be used to sample the equilibrium distribution
\eqref{eq:9}. In particular, the empirical distribution of the
trajectory, defined blow, 
\begin{equation}
\label{eq:empiricalmeas}
\lambda_T^\nu \equiv \frac 1T \int_B^{B+T} \delta_{(\bd x^\nu (t), \beta^\nu (t))} \ud t 
\end{equation}
will converge to the equilibrium distribution $\varrho$ as $T \to \infty$.
Note here $\delta_{(\bd x, \beta)}$ denotes the Dirac delta function,
that is the point mass function at $(\bd x, \beta)$, and $B$ is any
fixed burn-in time (which for simplicity will be taken as $B=0$ in the
numerical experiments).

To discuss the choice of $\nu$, it is important then to quantify the
speed of convergence of the empirical distribution
\eqref{eq:empiricalmeas} to the equilibrium: A better $\nu$ will
correspond to faster convergence and thus less simulation length of
the trajectory.  To quantify this convergence for STMD, the
theoretical tool we will use is the large deviation principle for the
empirical measure of stochastic processes. We will discuss the idea
and conclusion of the large deviation principle below, while defer the
rigorous definition and derivation to the Appendix.

When $T$ becomes large, the empirical distribution is expected to be
very close to the equilibrium. The large deviation principle
quantifies the probability that the empirical distribution is still
far away from the equilibrium: more specifically, let the probability
of the empirical distribution being equal to $\mu$ is on the order of
$O(e^{-I^{\nu}(\mu) T})$, with a rate functional, specific form given
below, $I^{\nu}(\mu) \geq 0$, and only vanish when $\mu$ is the
equilibrium distribution $\varrho$. This in particular tells us that
as $T \to \infty$, the likelihood that empirical distribution is
deviated away from the equilibrium is exponentially small, for which
the functional $I^{\nu}$ quantifies the rate of the exponential decay.
Hence a larger rate function indicates faster rate of convergence.

To specify the rate functional, let $\mu$ be a probability measure on
$\RR^{3n}\times\{\beta_1,\beta_2\}$ with smooth density and define
$\theta(\bd{x},\beta):=[d\mu/d\varrho](\bd{x},\beta)$, the ratio of
the probability density of $\mu$ and the equilibrium distribution. For
the simulated tempering process with switching frequency $\nu$, the
large deviation rate function for the empirical measure converging to
the stationary one is given by
\begin{equation}\label{eq:largedeviation}
I^\nu(\mu) = J_0(\mu)+\nu J_1(\mu),
\end{equation}
where
\begin{align}
  J_0(\mu)&=\sum_\beta \int \dfrac{1}{4\theta(\bd{x},\beta)^2} \left[ \beta^{-1}| \nabla_x \theta(\bd{x},\beta)|^2 \right] \mu(\ud\bd{x},\beta);\\
  J_1(\mu)&=\tfrac12 \sum_\beta \int g_{\beta\beta'}(\bd{x}) \biggl[ 1- \sqrt{\dfrac{\theta(\bd{x},\beta')}{\theta(\bd{x},\beta)}} \biggr]^2 \mu(\ud\bd{x},\beta).
\end{align}
Note that $J_1$ is non-negative and hence the large deviation rate
functional is a pointwise monotonic function in $\nu$. Thus, we
conclude from that a larger swapping rate $\nu$ corresponds to a
faster convergence of the empirical distribution to equilibrium.

\section{Infinite switching limit}\label{sec:inf}

As shown in the last section, the higher the switching rate $\nu$ is,
the faster the convergence of the empirical distribution, and thus the
sampling of the simulated tempering overdamped dynamics. On the other
hand however, a large $\nu$ value requires one to make many swapping
attempts, which slows down the actual simulation. To resolve this
problem, the key is that the limit when $\nu \to \infty$ can be taken
explicitly. This observation is first made for replica exchange
dynamics in \cite{dupuis2012infinite}. As for simulated tempering,
when the switching frequency between the temperatures $\nu \to \infty$,
the dynamics \eqref{eq:overdampmixt} converges to the following
\begin{equation}\label{eq:infswap1}
\dot{\bd{x}} = \bd{f}(\bd{x}) + \sqrt{2 (\beta_1^{-1} \omega_1(\bd{x}) + \beta_2^{-1} \omega_2(\bd{x}))} \bd{\eta}
\end{equation}
where the effective temperature
$T(\bd{x}) = \beta_1^{-1} \omega_1(\bd{x}) + \beta_2^{-1}
\omega_2(\bd{x})$
is a weighted average of $T_1 = \beta_1^{-1}$ and $T_2 = \beta_2^{-1}$, where the weight as a
function of $\bd{x}$ is given by \eqref{eq:weight}. Intuitively, this
can be understood as with fixed configuration $\bd{x}$ the proportion
of time the trajectory spends at temperature $\beta_j$ is given by
$\omega_j(\bd{x})$, and thus the effective temperature is given by a weighted average of the two temperatures. 

The Fokker-Planck equation for the probability density of $\bd{x}(t)$
corresponding to \eqref{eq:infswap1} can be written as
\begin{equation}\label{eq:fokkerplanck}
  \partial_t \rho(t) = \divop\bigl( \mathbb{B} (\rho(t) \grad U + k_B T_1 \grad \rho(t)) \bigr),
\end{equation}
where we have defined the effective potential
\begin{equation}\label{eq:effpotential}
  U(\bd{x}) = - \beta_1^{-1} \ln \varrho(\bd{x})
\end{equation}
and the mobility 
\begin{equation}
\mathbb{B}(\bd{x}) = \omega_1(\bd{x}) + \beta_1^{-1} \beta_2 \omega_2(\bd{x}).
\end{equation}
It is easy to check that the stationary solution of
\eqref{eq:fokkerplanck} is $\exp\bigl(-\beta_1 U(\bd{x}))$, which is
just $\varrho(\bd{x})$ by definition of $U$.

Note that this point of view also leads to a further simplification of
\eqref{eq:infswap1}, in the spirit of Lu and Vanden-Eijnden
\cite{lu2013infinite}. We may replace $\mathbb{B}$ in the
Fokker-Planck equation \eqref{eq:fokkerplanck} by a constant function,
which, for convenience, we will simply take to be the identity. This
substitution does not affect the stationary distribution, and hence
preserves the sampling property, but it changes the overdamped
Langevin dynamics associated to the Fokker-Planck equation.  It is
easy to see that this new dynamics is given by
\begin{equation}\label{eq:infswap2}
  \dot{\bd{x}} = \bigl(\omega_1(\bd{x}) + \beta_2 \beta_1^{-1} \omega_2(\bd{x})\bigr) \bd{f}(\bd{x}) + \sqrt{2 \beta_1^{-1}} \bd{\eta}.
\end{equation}
Compare with \eqref{eq:infswap1}, the noise term is additive in
\eqref{eq:infswap2}, and it still samples the desired stationary
distribution $\varrho$. Note that \eqref{eq:infswap2} is very similar
to the original overdamped Langevin dynamics \eqref{eq:overdamp}. The
only change is a scaling factor in front of the forcing term, which
involves the auxiliary (inverse) temperature $\beta_2$ and also the
weights $\omega_j$, $j = 1, 2$. Note that in practice
\eqref{eq:infswap2} can be implemented as an easy patch of existing
codes for molecular dynamics. 
 
\section{Generalizations}\label{sec:gen}

First, the idea in the previous section can be extended naturally to
more than two temperatures. For which, the dynamics in the infinite
switching limit is given by
\begin{equation}\label{eq:geninfswap}
\dot{\bd{x}} = \hbox{$\beta_1^{-1} \sum_{k} (\beta_k \omega_k(\bd{x})) \bd{f}(\bd{x}) + \sqrt{2 \beta_1^{-1}} \bd{\eta}$},
\end{equation}
where we define the weight here as
\begin{equation}\label{eq:genweight}
\begin{aligned}
\omega_k(\bd{x}) = \frac{n_k e^{-\beta_k V(\bd{x})}} {\sum_j n_j e^{-\beta_j V(\bd{x})}}, \quad k = 1, \cdots, N, 
\end{aligned}
\end{equation}
if $N$ is the number of temperatures used. Similarly, the dynamics \eqref{eq:geninfswap} can be
viewed as the overdamped dynamics with the effective potential given by \eqref{eq:effpotential} where the corresponding marginal equilibrium density becomes
\begin{equation}
\label{eq:Neqpdfmixt}
\varrho(\bd{x})　=　\frac{\sum_j n_j e^{-\beta_j V(\bd{x})}}{\sum_j n_j Z_{\beta_j}},
\end{equation}
as we have 
\begin{equation}
  \begin{aligned}
    \nabla_{\bd{x}} U(\bd{x}) &= - \beta_1^{-1} \frac{ \nabla_{\bd{x}} \sum_k n_j e^{-\beta_j V(\bd{x})} }{\sum_k n_j e^{-\beta_j V(\bd{x})} } \\
    & = \beta_1^{-1} \frac{ \sum_k \beta_k n_j e^{-\beta_j V(\bd{x})} \nabla V(\bd{x})}{\sum_k n_j e^{-\beta_j V(\bd{x})}} \\
    & = \beta_1^{-1} \sum_k \beta_k \omega_k(\bd{x}) \bd{f}(\bd{x}). 
  \end{aligned}
\end{equation}
Thus the gradient of the effective potential $U(\bd{x})$ is exactly
the forcing term in \eqref{eq:geninfswap}.

The infinite switching simulated tempering can be also generalized to
(underdamped) Langevin equation, rather than the overdamped Langevin
dynamics \eqref{eq:overdamp}; other thermostats can be also
used. Recall the Langevin equation 
\begin{equation}\label{eq:langevin}
\begin{cases}
\dot{\bd{x}} = m^{-1}\bd{p}, \\[2mm]
\dot{\bd{p}} = \bd{f}(\bd{x})-\gamma \bd{p}+\sqrt{2\gamma m\beta^{-1}}\bd{\eta},
\end{cases}
\end{equation}
where $m$ denotes the mass and $\gamma$ the friction coefficient, in
which case the generalization of \eqref{eq:infswap2} reads
\begin{equation}\label{eq:generalized}
\begin{cases}
  \dot{\bd{x}} = m^{-1}\bd{p}, \\[2mm]
  \dot{\bd{p}} = \beta_1^{-1} \sum_{k} (\beta_k \omega_k(\bd{x})) \bd{f}(\bd{x}) \\[1.5mm]
  \qquad \qquad -\gamma \bd{p}+\sqrt{2\gamma m\beta^{-1}}\bd{\eta}.
\end{cases}
\end{equation}
The structure of the equation is rather similar to the overdamped case
\eqref{eq:infswap2}, the only modification compared to
\eqref{eq:langevin} is the scaling factor in front of the forcing term
that amounts to a weighted average of the inverse temperatures.

\section{Connection with integrated tempering enhanced sampling}\label{sec:ITS}

The infinite switching limit of the simulated tempering sampling
scheme is very closely related to the integrated tempering enhanced
sampling (ITS) algorithm originally proposed
in~\cite{gao2008integrate}. The ITS algorithm introduces a temperature biased effective potential energy as  
\begin{equation}
V_{\text{eff}}(\bd{x}) = -\beta_1^{-1} \text{ln}\left( \sum_k n_{\beta_k} e^{-\beta_k V(\bd{x})} \right)
\end{equation}
and run molecular dynamics simulation on the surface. Here
$n_{\beta_k}$ is chosen to be some weighting factor for the
temperature $\beta_k$. Note that this is identically the same as we have in \eqref{eq:effpotential} with the marginal equilibrium $\varrho(\bd{x})$ given by \eqref{eq:Neqpdfmixt}, despite a constant difference which does not matter in sense of potential energy.

As a conclusion, the ITS algorithm can be viewed as the infinite
switching limit of the simulated tempering algorithm. As we discussed
in Section~\ref{sec:LDT}, the sampling efficiency of the
simulated tempering method increases as $\nu \to \infty$. Thus as a
corollary, the ITS algorithm is more efficient in sampling compared
with the simulated tempering algorithm at a finite switching rate. This will be further demonstrated in our numerical examples in the next section.



\section{Numerical examples}\label{sec:numerical} \label{sec:numerics}

\subsection{Simple high dimensional example}
We first consider a system in $D$ dimension moving on the following
potential with $\bd{x} = (x_0, \cdots, x_{D-1})$
\begin{equation}\label{eq:potentialD}
  V(\bd{x}) = (1 - x_0^2) ^2 - \dfrac14 x_0 + \sum_{j=1}^{D-1}\dfrac12 \lambda_j x_j^2
\end{equation}
where$\lambda_1, \lambda_2, \cdots, \lambda_{D-1}$ are parameters
controlling the stiffness of the harmonic potential in the
$x_0, x_1, \cdots, x_{D-1}$ directions.  We take the low (physical)
temperature $\beta_0 = 25$ and five artificial temperatures
$\beta_k = 25 \times 2^{-k}, k = 1,2,3,4,5$. For each $\beta_k$, its
weighting parameter $n_k$ are set to be inverse of the corresponding
partition function, i.e. $Z_{\beta_k}^{-1}$, which could vary when
number of dimensions changes. When simulating with STMD, we use
\begin{equation}\label{eq:st2}
\dot{\bd{x}} =  \beta_1^{-1} \beta(t) \bd{f}(\bd{x}) + \sqrt{2 \beta_1^{-1}} \bd{\eta}.
\end{equation}
so that its infinite swapping limit coincides with
\eqref{eq:geninfswap} (rather than a multi-temperature version of
\eqref{eq:infswap1} for the original overdamped dynamics
\eqref{eq:overdampmixt}). At finite switching rate, we just consider
the switching attempts of $\beta(t)$ from some $\beta_k$ to its
adjacent $\beta_{k-1}$ (if $k > 0$), or $\beta_{k+1}$ (if $k < 5$).
And in each attempt to switch, we shall first decide whether to switch
up or down with equal probability, if $0 < k < 5$.  The total
simulation time is $T_{tot} = 2.5 \times 10^6$ with time step
$dt = 0.025$, which means the total number of steps is
$N_{tot}= 10^8$. As mentioned previously, we will take no burn-in
period when estimating the average of physical observable.

\begin{figure}[htbp]
	\centering
	\includegraphics[width=210px]{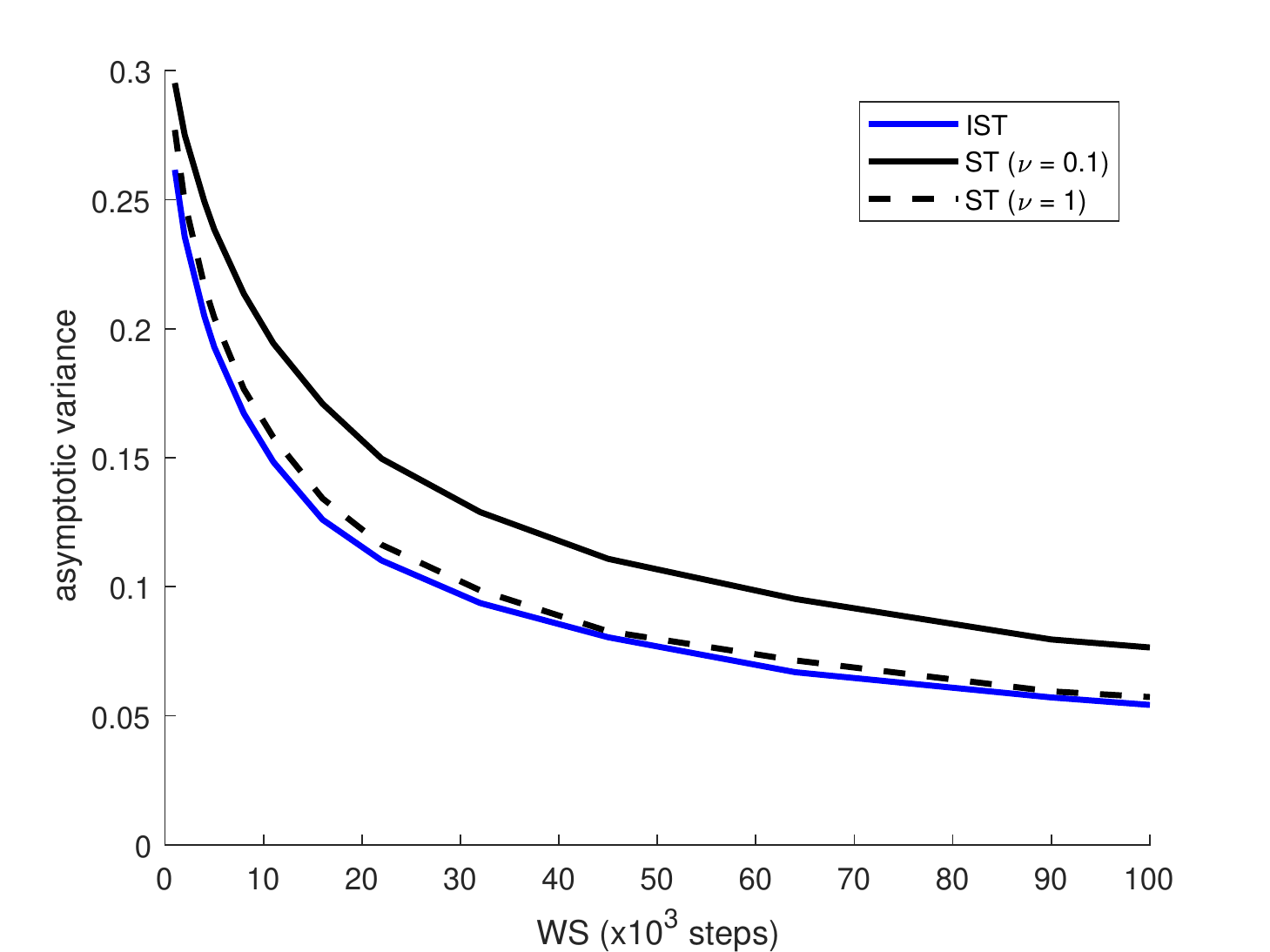}
	\caption{Asymptotic variance of $V(\bd{x})$ as in \eqref{eq:potentialD} ($D = 1$, thus a double well potential) for simulated tempering algorithms at different switching frequency and the infinite switching limit.}\label{fig:6T1DAV}
\end{figure}

\begin{figure}[htbp]
	\centering
	\includegraphics[width=210px]{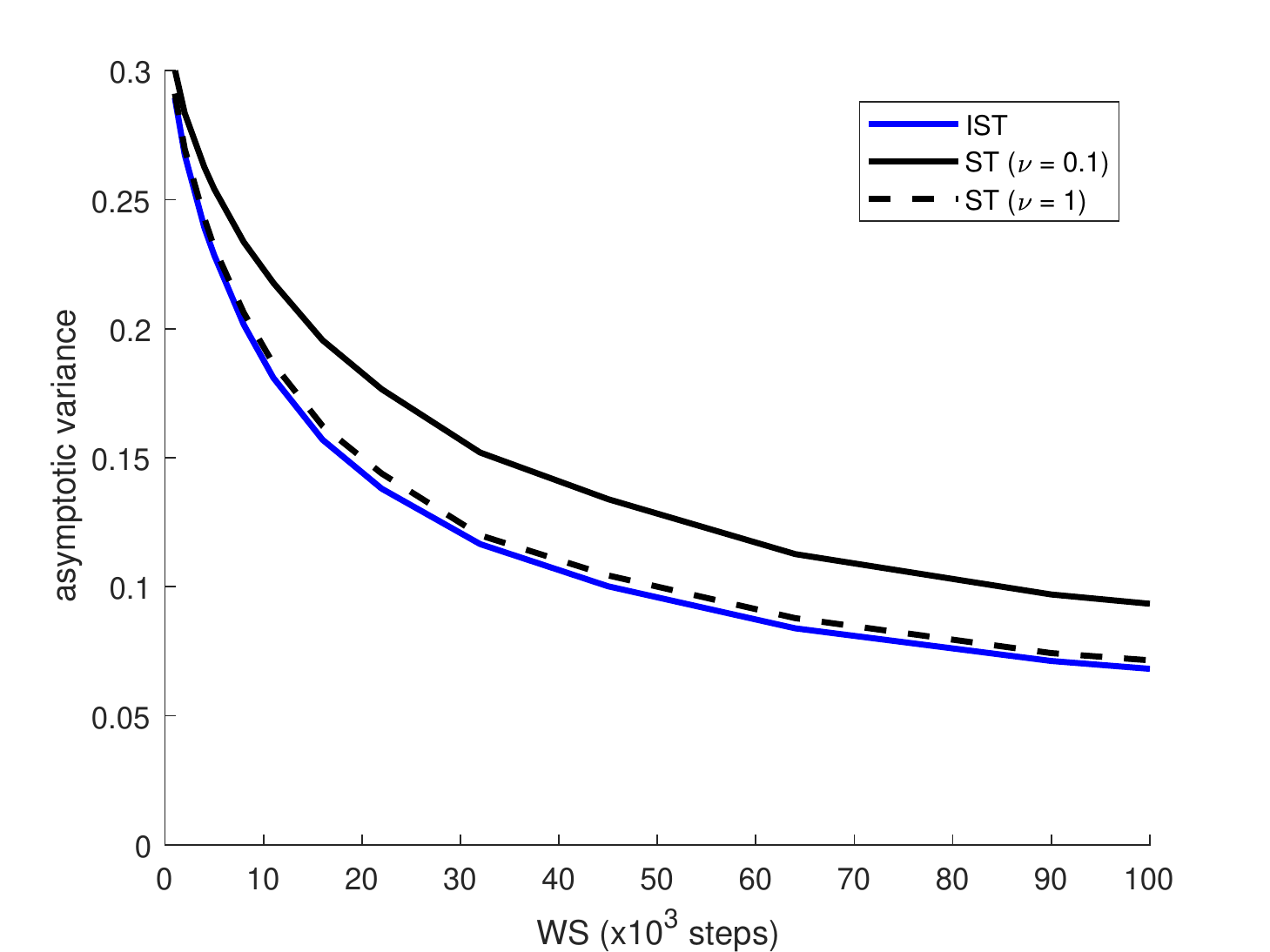}
	\caption{Asymptotic variance of $V(\bd{x})$ as in \eqref{eq:potentialD} ($D = 10$) for simulated tempering algorithms at different switching frequency and the infinite switching limit.}\label{fig:6T10DAV}
\end{figure}

To compare the algorithms, we calculate the asymptotic variance of the
observable $V(\bd{x})$  using a batch estimation
\begin{equation}\label{eq:asymptoticvariance}
AV = \text{Var}\;
\biggl\{ \sum_{t = (j-1)*WS + 1}^{j*WS}V(\bd{x}(t)), \; j = 1, \ldots, \frac{N_{\text{tot}}}{WS} \biggr\},
\end{equation}
where $AV$ stands for asymptotic variance, $WS$ stands for window size
of the batch. The results are shown in \figref{fig:6T1DAV} and
\figref{fig:6T10DAV} for $D = 1$ and $D = 10$ respectively. We observe
that for the finite switching rate, STMD at frequency $\nu = 1$ has a
lower asymptotic variance compared to $\nu = 0.1$, and moreover the
infinite swapping limit has a even lower asymptotic variance. We also
observe that for this example, the simulated tempering with $\nu = 1$
is already quite close to the infinite switching limit.

\subsection{Dimer in solvent example}

Now, to test the performance of our algorithm on a more realistic example, we apply \eqref{eq:generalized} for a dimer in solvent model as considered in 
\cite{dellago1999calculation}. This system consists of $N$ two-dimensional particles in a periodic box with side length $l$. All particles have the same mass $m$, and they interact with each other with the Weeks-Chandler-Anderson potential defined as
\begin{equation}\label{eq:WCApotential}
V_\text{WCA}(r)=4\epsilon \left( (\sigma /r)^{12} - (\sigma /r)^{6} \right) + \epsilon,
\end{equation}
if $r \le r_\text{WCA} = 2^{1/6}\sigma$, and $V_\text{WCA}(r) = 0$ otherwise, except for a pair of particles which interact via a double well potential
\begin{equation}\label{eq:dWpotential}
V_{dW}(r)=h \left( 1-\dfrac{(r-r_\text{WCA}-\omega)^2}{\omega^2} \right)^2.
\end{equation}
We take $N = 16$, $l = 4.4$, $\sigma = 1$, $h = 1$, $\omega = 0.5$, and $\epsilon = 1$ in the simulation. The physical temperature is $T_0 = 0.2$, and one artificial temperatures $T_1 = 1.0$.  The total simulation time is set to be $T_{tot} = 1.0 \times 10^5$ with time step $dt = 0.001$, which means the total number of steps is $N_{tot} = 10^8$. The quantity of interest is the free energy associated with the distance of the pair of particles interacting via the double well potential.

In this numerical example, we would still choose the weighting
parameter to be inverse partition function, which is however not
explicitly known or easily obtained a priori in this case. Therefore,
we use the following method to get approximation of these partition
functions, similar to that of \cite{gao2008integrate}, and then use
the estimates in the weighting factors $n_k$. First, a set of initial
guess $Z_k^{(0)}$ is chosen and used to run the dynamics
simulation. Then from a relatively short trajectory, we can add up the
corresponding weight to get the `proportion' of each temperature by
normalizing their summation to one. Since this proportion would go to
$1/2$ if the guess is accurate and the trajectory is infinitely long,
we can adjust our guess accordingly. More specifically, if the
proportions are $w_k,\ k = 0, 1$, the weighting factor will be updated
as
\begin{equation}\label{eq:newguess}
  Z_k^{\text{new}} =
\begin{cases}
Z_k^{\text{old}} \times 2 w_k, & \text{if } 2 w_k \in I; \\
Z_k^{\text{old}} \times \sqrt{2 w_k}, & \text{otherwise.}
\end{cases}
\end{equation}
where $I$ is a small neighborhood around $1$ to take into account of
the fluctuation due to the short simulation trajectory when estimating
$w_k$.  We iteratively update the partition function estimate
$\{Z_k^{(l)}\},\, l = 1, 2, \cdots$ until the proportions $w_k$ are
satisfactorily close to $1/2$.  Here we take $I = [0.35, 1.5]$,
initial guess $(Z_0^{(0)}, Z_1^{(0)}) = (1, 10^8)$ and number of
iterations $l_{\max} = 10$ with trajectory length $10^7$ steps each.

\begin{figure}[htbp]
	\centering
	\includegraphics[width=210px]{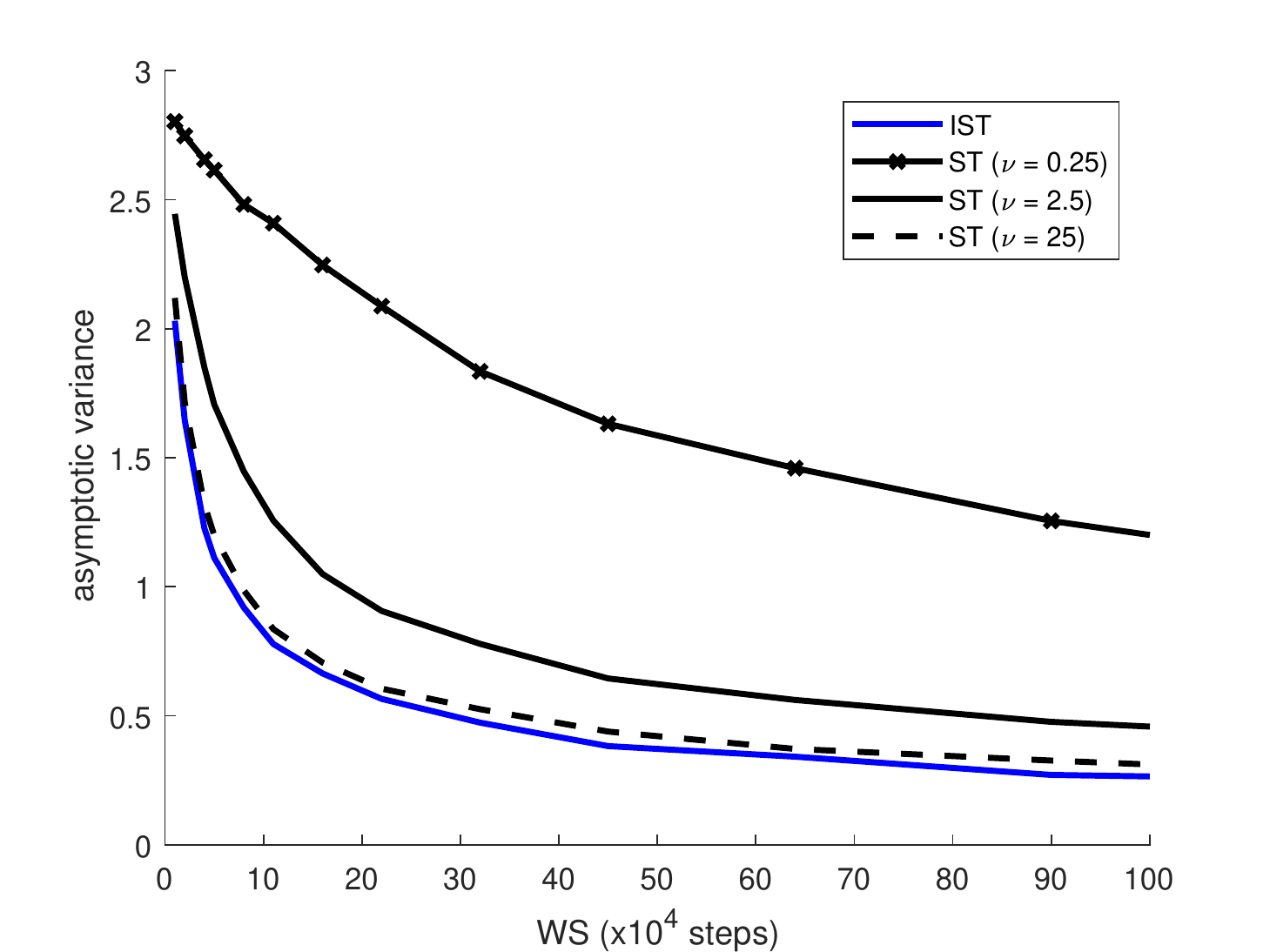}\\
	\caption{Asymptotic variance of the potential energy $V(\bd{x})$
      for the WCA example for simulated tempering algorithms at
      different switching frequency and the infinite switching
      limit.}\label{fig:WCAAV}
\end{figure}

The performance of STMD with various switching rate is shown \figref{fig:WCAAV}. It can be seen that the 
performance of STMD is better when switching rate increases from
$0.25$ to $25$; and the asymptotic variance converges to that of the IST as $\nu \to \infty$. This provides strong numerical validation of our theoretical results. 

\section{Conclusion}\label{sec:conclusion} \label{sec:conclude}

We justify that the sampling efficiency of the simulated tempering
method increases with the switching rate of the temperature, using
both the theoretical analysis based on large deviation of empirical
distribution and also numerical tests on two examples. This motivates
taking the infinite switching limit of the simulated tempering
dynamics, which recovers the integrated tempering enhanced sampling
method under a natural reformulation. The limiting dynamics can be
implemented as a patch to standard molecular dynamics by adding a
scaling factor to the force term based on a weighted average of all
temperatures involved. This leads to a practical scheme with higher
sampling efficiency than standard simulated tempering. 

\section*{Acknowledgment} 

HG is supported in part by National Science Foundation of China via grant 11622101. JL is supported in part by National Science Foundation via
grant DMS-1454939. JL and HG would also like to thank Yiqin Gao, Zhiqiang Tan, Eric
Vanden-Eijnden and Zhennan Zhou for helpful discussions.

\appendix

\section{Derivation of the large deviation rate functional
  $I^\nu(\mu)$}

To introduce formally the large deviation principle, let us introduce
some definitions. Let $S$ be a Polish space and $\cP (S)$ the space
of probability measures on $S$, equiped with the topology of weak
convergence. Under this weak topology, $\cP (S)$ itself is a Polish
space. Note that the empirical measure $\lambda_T^\nu$ is a
$\cP (S)$-valued random variable.

The large deviation principle for empirical distribution can be stated
as follows \cite{DonskerVaradhan:1975, DeuschelStroock:1989, DemboZeitouni}: 
A sequence of random probability measures $\{\gamma_T\}$ is said to satisfy a large deviation principle with rate function $I : \cP (S) \to [0, \infty]$, if for all open sets $O \subset \cP (S)$
\begin{equation}\label{eq:LDPopen}
\liminf_{T\to\infty} \frac 1T \log \PP (\gamma_T \in O) \ge - \inf_{\mu \in O} I (\mu),
\end{equation}
for all closed sets $C \subset \cP (S)$
\begin{equation}\label{eq:LDPclosed}
\liminf_{T\to\infty} \frac 1T \log \PP (\gamma_T \in C) \le - \inf_{\mu \in C} I (\mu),
\end{equation}
and for all $M < \infty$, $\{\mu: I (\mu) \le M\}$ is compact in
$\cP (S)$.  In particular, once we obtain the rate functional,
\eqref{eq:LDPopen} and \eqref{eq:LDPclosed} quantifies how unlikely
the empirical distribution is far from the equilibrium when $T$ is
large, which was used in this work to justify the infinite switching
limit.

To find the rate functional $I^\nu$ for the simulate tempering
overdamped dynamics with switching rate $\nu$, we divide the
infinitesimal generator into two parts
\begin{equation}\label{eq:infgendiv}
(\mathcal{L}^\nu u)(\bd{x},\beta)= (\mathcal{L}^\nu_\mathrm{diff} u + \nu \mathcal{L}^\nu_\mathrm{jump} u)(\bd{x},\beta)\\
\end{equation}
where
\begin{equation}
\begin{cases} 
\mathcal{L}^\nu_\mathrm{diff} = -\nabla_{\bd{x}} V(\bd{x})\cdot\nabla_{\bd{x}} + \beta^{-1}\nabla_{\bd{x}} \cdot\nabla_{\bd{x}},\\
\mathcal{L}^\nu_\mathrm{jump} = -g_{\beta,\beta'}(\bd{x}) + g_{\beta',\beta}(\bd{x}).
\end{cases}
\end{equation}
The rate functional thus has an additive structure
\eqref{eq:largedeviation} with the rate functional $J_0$ and $J_1$
correspond to $\mathcal{L}^\nu_\mathrm{diff}$ and
$\mathcal{L}^\nu_\mathrm{jump}$, respectively. 

To get an explicit formula of the rate functionals, we first consider
$J_0$. Let $ f(\bd{x},\beta) = \sqrt{\theta(\bd{x},\beta)}$, then following 
Donsker and Varadan \cite{DonskerVaradhan:1975}, we have
\begin{equation}\label{eq:derivationJ0}
\begin{aligned}
J_0(\mu) & = \sum_{\beta} \int f (\mathcal{L}_{\mathrm{diff}}^{\nu} f) \,\varrho(\ud\bd{x}, \beta)\\
 = & -\sum_{\beta} \int \biggl[ -(\nabla_{\bd{x}} V(\bd{x})\cdot \nabla_{\bd{x}} f) f \\
  & \qquad \qquad \qquad + (\beta^{-1}\nabla_{\bd{x}} \cdot\nabla_{\bd{x}} f) f \biggr] \varrho(\ud\bd{x},\beta)\\
=& -\sum_{\beta} \frac 12 Z_\beta^{-1} \left\{\int -(\nabla_{\bd{x}} V(\bd{x})\cdot \nabla_{\bd{x}} f) f e^{-\beta V(\bd{x})} \ud\bd{x}\right. \\
&\qquad \qquad \qquad + \left. \beta^{-1}\int(\nabla_{\bd{x}} \cdot\nabla_{\bd{x}} f) f e^{-\beta V(\bd{x})} \ud\bd{x}\right\}\\
=& -\sum_{\beta} \frac 12 Z_\beta^{-1} \left\{\int -(\nabla_{\bd{x}} V(\bd{x})\cdot \nabla_{\bd{x}} f) f e^{-\beta V(\bd{x})} \ud\bd{x}\right. \\
&\qquad \qquad \qquad- \left. \beta^{-1}\int\nabla_{\bd{x}} f \cdot \nabla_{\bd{x}} ( f e^{-\beta V(\bd{x})}) \ud\bd{x}\right\}\\
=& -\sum_{\beta} \frac 12 Z_\beta^{-1} \left\{-\beta^{-1}\int \abs{\nabla_{\bd{x}} f}^2 e^{-\beta V(\bd{x})} \ud\bd{x}\right\}\\
=& \sum_{\beta} \int \frac 12\beta^{-1}\abs{\nabla_{\bd{x}} \sqrt{\theta(\bd{x},\beta)}}^2 \rho_{\beta}(\ud\bd{x})\\
=& \sum_\beta \int \dfrac{1}{8\theta(\bd{x},\beta)} \left[ \beta^{-1}| \nabla_x \theta(\bd{x},\beta)|^2 \right] \rho_{\beta}(\ud\bd{x}).
\end{aligned}
\end{equation}

Next we consider the rate functional $J_1$ corresponding to the jump
process of the temperature with generator
$\mathcal{L}_{\mathrm{jump}}^{\nu}$. From the definition of
$\theta(\bd{x},\beta)$, we have
\begin{equation}\label{eq:lemma0}
\theta(\bd{x}, \beta) = \dfrac{2\mu(\bd{x},\beta)}{\rho_\beta(\bd{x})},
\end{equation}
and therefore
\begin{equation}
  \begin{aligned}
\dfrac{\mu(\bd{x},\beta')}{\mu(\bd{x},\beta)} &= \dfrac{\theta(\bd{x}, \beta')}{\theta(\bd{x}, \beta)} \cdot \dfrac{\rho_{\beta'}(\bd{x})}{\rho_\beta(\bd{x})}\\ 
& =  \dfrac{\theta(\bd{x}, \beta')}{\theta(\bd{x}, \beta)} \cdot \dfrac{g_{\beta,\beta'}(\bd{x})}{g_{\beta',\beta}(\bd{x})}.
\end{aligned}
\end{equation}
Combined with the Cauchy-Schwartz inequality, we obtain that for any choice of non-negative $u(\bd{x}, \beta)$
\begin{equation}\label{eq:lemma1J1}
\begin{aligned}
  & \int \bigg( g_{\beta',\beta}(\bd{x}) \dfrac{u(\bd{x}, \beta')}{u(\bd{x}, \beta)} \mu(\ud\bd{x}, \beta)\\
  & \qquad + g_{\beta,\beta'}(\bd{x}) \dfrac{u(\bd{x}, \beta)}{u(\bd{x}, \beta')} \mu(\ud\bd{x}, \beta') \bigg)\\
  & \ge \int 2 \sqrt{g_{\beta',\beta}(\bd{x})g_{\beta,\beta'}(\bd{x})\dfrac{\mu(\bd{x}, \beta')}{\mu(\bd{x}, \beta)}} \mu(\ud\bd{x}, \beta)\\
  & = \int 2 g_{\beta,\beta'}(\bd{x}) \sqrt{
    \dfrac{\theta(\bd{x},\beta')}{\theta(\bd{x},\beta)} }
  \mu(\ud\bd{x}, \beta),
\end{aligned}
\end{equation}
where the equality is attained if and only if for each
$\bd{x} \in \RR^{3n}$ and each $\beta$ it holds
\begin{equation}\label{eq:lemma2J1}
\dfrac{u(\bd{x}, \beta')}{u(\bd{x}, \beta')} \propto \dfrac{g_{\beta,\beta'}(\bd{x})\mu(\bd{x}, \beta')}{g_{\beta',\beta}(\bd{x})\mu(\bd{x}, \beta)}
= \dfrac{\rho_{\beta'}(\bd{x})\mu(\bd{x}, \beta')}{\rho_{\beta}(\bd{x})\mu(\bd{x}, \beta)}.
\end{equation}
In particular, this gives the choice of $u(\bd{x}, \beta)$ to make 
equality holds in \eqref{eq:lemma1J1}.

Using the variational characterization of the large deviation rate
functional as in \cite{DonskerVaradhan:1975}, we have
\begin{equation}\label{eq:derivationJ1}
\begin{aligned}
 J_1(\mu) & =  -\inf_{u\geq 0} \sum_{\beta} \int \left( \dfrac{\mathcal{L}^\nu_\mathrm{jump} u}{u} \right) (\bd{x}, \beta) \, \mu(\ud\bd{x}, \beta)\\
=   & - \inf_{u\geq0} \sum_{\beta} \int \biggl[ -g_{\beta,\beta'}(\bd{x}) \\
  & \qquad \qquad \qquad + g_{\beta',\beta}(\bd{x}) \dfrac{u(\bd{x}, \beta')}{u(\bd{x}, \beta)} \biggr] \mu(\ud\bd{x}, \beta)\\
=   & \frac12 \sum_{\beta} \int \left[ g_{\beta,\beta'}(x) \mu(\ud\bd{x}, \beta) + g_{\beta',\beta}(\bd{x}) \mu(\ud\bd{x}, \beta') \right]\\
& - \frac12 \inf_{u\geq 0} \sum_{\beta} \int \bigg[ g_{\beta',\beta}(\bd{x}) \dfrac{u(\bd{x}, \beta')}{u(\bd{x}, \beta)} \mu(\ud\bd{x}, \beta)\\
& \qquad \qquad + g_{\beta,\beta'}(\bd{x}) \dfrac{u(\bd{x}, \beta)}{u(\bd{x}, \beta')} \mu(\ud\bd{x}, \beta') \bigg]\\
\overset{\eqref{eq:lemma1J1}}{=}   & \frac12 \sum_{\beta} \int g_{\beta,\beta'}(x) \left[ 1 + \dfrac{\theta(\bd{x},\beta')}{\theta(\bd{x},\beta)} \right] \mu(\ud\bd{x}, \beta)\\
& - \frac12 \sum_{\beta} \int 2 g_{\beta,\beta'}(\bd{x}) \sqrt{ \dfrac{\theta(\bd{x},\beta')}{\theta(\bd{x},\beta)} } \mu(\ud\bd{x}, \beta)\\
=   & \frac12 \sum_\beta \int g_{\beta\beta'}(\bd{x}) \left[ 1- \sqrt{\dfrac{\theta(\bd{x},\beta')}{\theta(\bd{x},\beta)}} \right] ^2 \mu(\ud\bd{x},\beta).
\end{aligned}
\end{equation}
Therefore, we arrive at the explicit formula for the large deviation
rate functionals. We remark that a rigorous proof of the rate
functional derived from the above calculation can be found in
\cite{DupuisLiu}.

\bibliography{isst}
\end{document}